\documentclass[review, 11pt]{elsarticle}
\usepackage{hyperref}
\usepackage{amsmath, amsthm, amscd, amsfonts, amssymb, graphicx, color}
\usepackage{amsmath,amsthm,amscd,amsfonts,amssymb,enumerate,pgf,tikz,mathrsfs,caption}
\usepackage{tikz-cd}
\usepackage{caption}
\usepackage{hyperref}
\hypersetup{nesting=true,debug=true,naturalnames=true}
\usepackage{graphicx,amssymb,upref}
\usepackage[all]{xy}
\theoremstyle{plain}
\newtheorem{theorem}{Theorem}[section]
\newtheorem{lemma}[theorem]{Lemma}

\theoremstyle{definition}
\newtheorem{definition}[theorem]{Definition}
\newtheorem{example}[theorem]{Example}

\newtheorem{corollary}[theorem]{Corollary}

\theoremstyle{remark}

\numberwithin{equation}{section}

\begin{document}

\begin{frontmatter}

\title{Partial semigroup partial dynamical systems and Partial Central Sets }

\author{H. Goodarzi, M. A. Tootkaboni$^{*}$}
\cortext[mycorrespondingauthor]{M. A. tootkaboni}
\ead{tootkaboni.akbari@gmail.com}
\address{Department of Pure Mathematics, Faculty of Mathematical Sciences,\\
 University of Guilan,\\
  Rasht-Iran.}
\author{Arpita Ghosh}
\ead{arpi.arpi16@gmail.com}
%\address{Pure and Applied Science, Midnapore City College, Kuturia, Bhadutala,\\
% Paschim Medinipur,\\
%  West Bengal-721129, India\\
%Email: souravkantipatra@gmail.com}

\begin{abstract}
H. Furstenberg defined Central sets in $\mathbb{N}$ by using the notions of topological dynamics, later Bergelson and Hindman characterized central sets in $\mathbb{N}$ and also in  arbitrary  semigroup in terms of algebra of Stone-\v{C}ech compactification of that set. We state the new notion of large sets in a partial semigroup setting and characterize the algebraic structure of the sets by using the algebra of Stone-\v{C}ech compactification. By using these notions, we introduce the \emph{Partial Semigroup Partial Dynamical System(PSPDS)} and show that topological dynamical characterization of central sets in a partial semigroup is equivalent to the usual algebraic characterization.
\end{abstract}


\end{frontmatter}

%\linenumbers

\section{Introduction}
Central sets play an extensive role in Ramsey's theory. It is primarily introduced by Furstenberg in terms of proximal and uniformly recurrent points of a  topological dynamical system on natural numbers, $\mathbb{N}$ in \cite{Furstenberg}. After that, in \cite{Bergelson and Hindman}, Bergelson and Hindman explore an algebraic characterization of such sets using the  algebra of Stone-\v{C}ech compactification of  $ \mathbb{N}$,  and subsequently, they also generalize for any arbitrary discrete semigroups.

A topological dynamical system on a compact topological space $X$ is a set of continuous endomorphisms $\{T_s\}_{s\in S}$ on $X$, parametrized by some semigroup $(S, \ast)$, satisfying the rule $T_s \circ T_t = T_{s\ast t}$ for all $s, t \in S.$ Therefore, it is consist of a semigroup $S$ and an  compact space $X$. Given a discrete semigroup $S$, the associated  Stone-\v{C}ech compactification $\beta S$, the set of all ultrafilters on $S$, forms a compact topological space generated by the basis elements $\bar{A}:= \{p \in \beta S: A \in p\}$ for $A \subseteq S.$  Note that $S \subseteq \beta S$ identified as the set of all principal ultrafilters. The operation on $S$ extends to $\beta S$ in such a way that for each $s \in S$ and $u \in \beta S$ the self maps 
$\lambda_s : p \mapsto s\ast p$ and $ \lambda_u : p \mapsto u \ast p $ on $\beta S$ are continuous. Shi and Yang studied dynamical systems on the compact space of continuous functions from $S$ to $\{0,1\}$ both having discrete topologies and proved certain dynamical classifications of central sets in the same spirit of Furstenburg, see \cite{Shi and Yang}. To understand it precisely, let us recall that for a dynamical system $(X, \{T_s\}_{s\in S})$, a point $y \in X$ is called \emph{uniformly recurrent} if and only if for every neighborhood $U$ of $y$, there exists a finite subset $F$ of $S$ with 
\[
\bigcup_{f \in F} f^{-1}\{s \in S: T_s(y) \in U\}=S.
\]
Two points $x, y \in X$ are called \emph{proximal} if and only if for every neighborhood $U$ of $\{(w,w): w\in X\}$, there exists $s\in S$ such that $(x,y) \in (T_s \times T_s)^{-1}(U).$ With this definitions, Shi and Yang's result states that a subset $C$ of $S$ is central if and only if $C$ is \emph{dynamically central}, that is,  there exists a dynamical system $(X, \{T_s\}_{s\in S})$, elements $x, y \in X$, and a neighborhood $U$ of $y$ such that $y$ is uniformly recurrent, $x$ and $y$ are proximal and $C= \{s \in S: T_s(x) \in U\}$ in \cite{Shi and Yang}.

The notion of partial semigroups is another important notion in the study of Ramsey
theory. Central sets in partial semigroups first occur in the work of McLeod \cite{McLeod1}, and  has been extensively studied by several authors \cite{Ghosh}, \cite{Hindman and Pleasant}. In this article, we introduce two frameworks of topological dynamical systems for partial semigroups, say \emph{partial semigroup partial dynamical system} (PSPDS) and \emph{partial semigroup dynamical system} (PSDS), consequently, establish a few new notions of large sets in partial semigroups and characterize the algebraic structure of these sets by using the algebra of Stone-\v{C}ech compactification. Finally, we prove that algebraic central sets are equivalent to certain (partial) dynamically central sets in the setting of PSPDS. 

\emph{Partial semigroup dynamical system} (PSDS), has been recognized in research in two ways: as a monotone action, see Definition 19.30 in \cite{Hindman} and as an $IP$-system see \cite{Fur-1987}.

\section{Preliminary}
In this section, we state some essential concepts which is needed them in our work. The first, we concentrate on the concept of the adequate partial semigroup. For more details, see \cite{Hindman,McCutcheon}.
\begin{definition}(Partial Semigroup) 
	Let $S$ be a  non-empty set.
	A pair $(S, \bullet)$ is called partial semigroup, when $\bullet$ is an operation defined on a subset $D$ of $S \times S$ to $S$ and satisfies the statement
	that for all $x, y, z \in S$, $(x \bullet y) \bullet z = x \bullet (y \bullet z)$ in the sense that if either side is defined,
	so is the other, and they are equal. 	
\end{definition}
We say that $x \bullet y$ is defined if $(x, y) \in D$. For every $ x \in S$, we define $R(x)= R_S(x) =\{s \in S : x \bullet s ~\text{is defined in S}\}$ and,  
$L(x)=\{s\in S:s\bullet x\mbox{ is well defined}\}$. 

A subset $I$ of $S$ is called left ideal if $y \bullet x \in I$ for $y \in L(x)$ and $x \in I$. The right ideal is defined similarly. We say that $I$ is an ideal if it is both a left ideal and a right ideal. A subset $L$ of $S$ is a minimal left ideal if $L$ is a left ideal of $S$ and, whenever $J$ is a left ideal of $S$ and $J \subseteq L$ implies that $J=L$. The minimal right ideal defined similarly.  An element $p \in S$ is an idempotent if $p \bullet p = p$ and the collection of all idempotents is denoted by $E(S)$. 

\begin{example}{\label{2.2}}
	\begin{flushleft}
		(i) Let $\mathcal{P}_f(\mathbb{N})= \{F : \emptyset \ne F \subseteq \mathbb{N} ~ and~ F~ \text{is finite}\}$. Then $(\mathcal{P}_f(\mathbb{N}), \biguplus)$ is a partial semigroup, when $F\biguplus G= F\cup G$ if $F\cap G=\emptyset$.\\
		(ii) Let $\{x_n\}_{n=1}^{\infty}$ be a sequence in an arbitrary semigroup $(S, \cdot )$. Define 
		\[
		T= FP(\{x_n\}_{n=1}^{\infty}) = \{\prod_{n \in F} x_n : F \in \mathcal{P}_{f}(\mathbb{N})\}.
		\]
	\end{flushleft}
	Then $(T,\bullet)$ is a partial semigroup if $"\bullet "$ is defined by
	\begin{equation*}
	(\prod_{n \in F} x_n) \bullet (\prod_{n \in G} x_n) =  \left\{
	\begin{array}{rl}
	\prod_{n \in F \cup G} x_n & \text{if}~  max F <  min G\\
	\text{undefined} & \text{if}~	max F \geq  min G
	\end{array} \right.
	\end{equation*}
\end{example}

Now, we recall some of the basic properties of an adequate partial semigroup.
\begin{definition}
	Let $(S, \bullet)$ be a partial semigroup.
	\begin{flushleft}
		(a)
		For $H \in \mathcal{P}_f(S),$ we define $R(H) = \bigcap_{s \in H} R(s)$, and $L(H) = \bigcap_{s \in H} L(s)$\\
		(b)
		We say that $(S, \bullet)$ is right adequate (left adequate) if  $R(H)\ne \emptyset\,\,(L(H)\neq\emptyset)$ for all $H \in \mathcal{P}_f(S)$. We say that $(S,\bullet)$ is adequate if $S$ is right and left adequate partial semigroup. \cite{Bergelson5}\\
		(c)
		$\delta_R S = \bigcap_{H \in \mathcal{P}_f(S)} c\ell_{\beta S}(R(H))$ and $\delta_L S = \bigcap_{H \in \mathcal{P}_f(S)} c\ell_{\beta S}(L(H))$.
	\end{flushleft}
\end{definition}
If $(S, \bullet)$ is partial semigroup, so is $(\beta S, \bullet)$ (see Proposition B.3 of \cite{Krautzberger}). Now we state some facts about partial semigroups and right adequate partial semigroups. By a similar way, it is true for left adequate partial semigroups. In other words, $(S, \bullet)$ is a right adequate partial semigroup if the collection $\{R(x):x\in S\}$ has the finite intersection property.  A subset $\delta_R S $ of $\beta S$ is a compact right topological semigroup (by Theorem 2.10 of \cite{McCutcheon}). Therefore, every concept for a compact right topological semigroup will be true for $\delta_R S$ (For more detail, see \cite{McLeod1}). If $(S,\bullet)$ is a right adequate partial semigroup, for every $s\in S$ and $A\subseteq S$, we define 
\[
s^{-1}A=\{t \in R(s): s \bullet t \in A\}.
\]

\begin{lemma}
	Let $(S, \bullet)$ be a partial semigroup, let $A \subseteq S$ and let $a, b, c \in S.$ Then $c \in b^{-1}(a^{-1}A)$ if and only if both $b \in R(a)$ and $c \in (a \bullet b)^{-1}A.$ In particular, if $b \in R(a),$ then $b^{-1}(a^{-1}A)= (a \bullet b)^{-1}A.$
\end{lemma}

\begin{proof}
	See Lemma 2.3 in \cite{McCutcheon}. 
\end{proof}

\begin{definition}
	Let $(S, \bullet)$ be a right adequate partial semigroup.
	\begin{flushleft}
		(a)
		For $a \in S$ and $q \in \overline{R(a)},$ we define $a \bullet q=\{A \subseteq S \mid a^{-1}A \in q\}.$\\
		(b)
		For $p,q \in \beta S$, if $p\bullet q$ is well defined, then we define
		\[
		p \bullet q =\{A \subseteq S \mid \{a^{-1}A \in q\}\in p\}.
		\]
	\end{flushleft}
\end{definition}

\begin{lemma}
	Let $(S, \bullet)$ be a right adequate partial semigroup.
	\begin{flushleft}
		(a)
		If $a \in S$ and $q \in \overline{ R(a)},$ then $a \bullet q \in \beta S.$\\
		(b) 		
		Let $p \in \beta S,$ $q \in \delta_R S$ and $a \in S.$ Then $R(a) \in p \bullet q$ if and only if $R(a) \in p$.\\
		(c)
		If $p, q \in \delta_R S,$ then $p \bullet q \in \delta_R S$.\\
	\end{flushleft}
\end{lemma}
\begin{proof}
	See Lemma 2.7 in \cite{McCutcheon}.
\end{proof}
\begin{lemma}{\label{2.7}}
	Let $S$ be a partial semigroup. \\
	\begin{flushleft}
		(a) Let $L_1$ and $L_2$ be left ideals of $S$. Then $L_1\cap L_2$ is left ideal if and only if $L_1\cap L_2\neq\emptyset$.\\
		(b) Let $L$ be a left ideal and $R$ be a right ideal. If $(\bigcup_{s\in L}L(s))\cap R\neq\emptyset$ then $L\cap R\neq\emptyset$.\\
		(c) Let $s\in S$. Then, $s\bullet S$ is  a right ideal, $S\bullet s$ is a left ideal, and $S\bullet s\bullet S$ is an ideal.\\
		(d) Let $e\in E(S)$. Then $e$ is a left identity for $e\bullet S$, a right identity for $S\bullet e$, and an identity for $e\bullet S\bullet e$.
	\end{flushleft}
\end{lemma}
\begin{proof}
	It is obvious.
\end{proof}
\begin{definition}{\label{def-p.g}}
	We say that $(G,\bullet)$ is a partial group if
	\begin{flushleft}
		(a) $(G,\bullet)$ is a partial semigroup.\\
		(b) There exists a left identity $e\in G$, i.e., there exists $e\in E(G)$ such that for every $s\in G$, $e\in L(s)$ and $e\bullet s=s$, and\\
		(c) For every $x\in G$, there exists $y\in L(x)$ such that $y\bullet x=e$, ( $y$ is called left $e$-inverse for $x$) . 
	\end{flushleft}
\end{definition}
By a similar way, right identity and right $e$-inverse have been defined.
\begin{lemma}{\label{2.9}}
	Let $(S,\bullet)$ be a partial semigroup. The following statements are 
	equivalent:
	\begin{flushleft}
		(a) $(S,\bullet)$ is a partial group.\\
		(b) There is a two sided identity $e$ for $S$ with the property that for each $x\in S$ there is some $y\in R(x)\cap L(x)$ such that $y$ is a (two sided) $e$-inverse for $x.$\\
		(c) There is a left identity for $S$ and given any left identity $e$ for $S$ and any $x\in S$ there is some $y\in L(x)$ such that $y$ is a left $e$-inverse for $x$.
	\end{flushleft}
\end{lemma}
\begin{proof}
	(a) implies (b). Since $S$ is a partial group, let $e$ be the left identity of $S$. We know that for every $x\in S$, there exists $y\in L(x)$ such that $y\bullet x=e$. Now, let $z\in L(y)$ be a left $e$-inverse for $y$. Then, by the below relations, we have $y\in L(x)\cap R(x)$.
	\begin{align*}
	x\bullet y=&e\bullet (x\bullet y)\\
	=&(z\bullet y)\bullet(x\bullet y)\\
	=&z\bullet(y\bullet(x\bullet y))\\
	=&z\bullet((y\bullet x)\bullet y)\\
	=&z\bullet(e\bullet y)\\
	=&z\bullet y\\
	=&e,
	\end{align*}
	and so $y$ is a right $e$-inverse for $x$. 
	
	Now, let $x\in S$ be given. Pick $y\in L(x)\cap R(x)$ as $e$-inverse for $x$. So, by   
	\[
	x\bullet e=x\bullet(y\bullet x)=(x\bullet y)\bullet x=e\bullet x=x,
	\]
	we have $e\in L(x)$ and so $e$ is right identity for $x$.\\
	(b) implies (c). It is obvious.\\
	(c) implies (a). Trivially.
\end{proof}
\begin{lemma}
	Let $(S,\bullet)$ be a partial semigroup. The following statements are equivalent:
	\begin{flushleft}
		(a) $(S,\bullet)$ is a partial group.\\
		(b) There is a right identity $e$ for $S$ such that for each $x\in S$ there is some $y\in R(x)$ such that $y$ is a right $e$-inverse for $x$.\\
		(c) There is a right identity $e$ for $S$ and given any right identity $e$ for $S$ and any $x\in S$ there is some $y\in R(x)$ such that $y$ is a right $e$-inverse for $x$.
	\end{flushleft}
\end{lemma}
\begin{proof}
	By Lemma \ref{2.9}, it is obvious.
\end{proof}
\begin{theorem}{\label{1.40}}
	Let $S$ be a partial semigroup and let $e$ be a left identity for $S$ such that for each $x\in S$ there is some $y\in R(x)$ with $x\bullet y = e$. Let $Y=E(S)$ and let $G=S\bullet e$. Then $Y$ is a partial right zero semigroup, $G$ is a partial  group, and $S= G\bullet Y$ and $G\times Y$ are partially isomorphic.
\end{theorem}
\begin{proof}
	The proof is alike Theorem 1.40 in \cite{Hindman}. We have to state the proof.
	
	The first we prove $x\bullet y=y$ for every $x\in Y$ and every $y\in R(x)$. Let $x\in Y$ and $y\in R(x)$ be given. Now, pick $z\in R(x)$ such that $x\bullet z=e$. Thus, $
	x\bullet e=x\bullet x\bullet z=x\bullet z=e$. Therefore $x\bullet y=x\bullet(e\bullet y)=(x\bullet e)\bullet y=e\bullet y=y$, proved our claim. Therefore for every $x\in Y$ and for every $y\in L(x)\cap Y$, $x\bullet y=y\in Y$, and $Y\neq\emptyset$ because $e\in Y$. Hence, $Y$ is a partial right zero semigroup.
	
	Now, by Lemma \ref{2.7}, $e$ is a right identity for $G=S\bullet e$. Also, if $x\in G$, then there exists $y\in R(x)$ such that $x\bullet y=e$. Therefore $x\bullet y\bullet e=e\bullet e=e$ and $y\bullet e\in G$. Also, $G\bullet G=S\bullet e \bullet S\bullet e\subseteq S\bullet S\bullet S\bullet e\subseteq S\bullet e=G$, it follows that $G$ is a partial group.
	
	Now, define $\varphi:D\subseteq G\times Y\rightarrow S$ by $\varphi(g,y)=g\bullet y$, where 
	\[
	D=\{(g,y):g\bullet y\mbox{ is well defined}\}.
	\] 
	We show that $\varphi$ is a partial homomorphism. In order to, let $(g_1,y_1),(g_2,y_2)\in D$, then,
	\begin{align*}
	\varphi(g_1,y_1)\bullet\varphi(g_2,y_2)=&(g_1\bullet y_1)\bullet (g_2\bullet y_2)\\
	=&	g_1\bullet g_2\bullet y_2\\
	=&	g_1\bullet g_2\bullet y_1\bullet y_2\\
	=&\varphi(	g_1\bullet g_2, y_1\bullet y_2)
	\end{align*} 
	Now, we show that $\varphi$ is surjective. Let $s\in S$, so $s\bullet e\in S\bullet e=G$, and so there exists $x\in S\bullet e$ such that $x\bullet (s\bullet e)=(s\bullet e)\bullet x=e$. We prove that $x\bullet s\in Y=E(S)$. Since, 	 
	
	\begin{align*}
	x\bullet s\bullet x\bullet s=&x\bullet s\bullet e\bullet x\bullet s (\mbox{because }x\in G,\,\, e\bullet x=x)\\
	=&x\bullet e\bullet s\\
	=&x\bullet s(\mbox{ because }x\in G,\,\, x\bullet e=x).
	\end{align*}
	Therefore $(s\bullet e,x\bullet s)\in D$ and $\varphi(s\bullet e,x\bullet s)=s\bullet e\bullet x\bullet s=e\bullet s=s$. Since $\varphi$ is onto, we established that $S=G\bullet Y$.
	
	Now, we prove that $\varphi$ is one-to-one. Let $(g,y)\in D$, and let $s=\varphi(g,y)$. We show that $g=s\bullet e$ and $y=x\bullet s$, where $x$ is the unique inverse of $s\bullet e$ in $G$. Since $s=g\bullet y$ so
	\begin{align*}
	s\bullet e=&g\bullet y\bullet e\\
	=&g \bullet e \,\,\,\,\,(y\bullet e=e)\\
	=&g\,\,\,\,\,(\mbox{ because }g\in G).
	\end{align*}
	Also
	\begin{align*}
	x\bullet s=&x\bullet g\bullet y\\
	=&x\bullet g\bullet y\bullet e\bullet y\,\,\,\,\,(y\bullet e\in Y\mbox{ so } y\bullet e\bullet y=y)\\
	=&x\bullet s\bullet e\bullet y\\
	=&e\bullet y\\
	=& y.
	\end{align*}
\end{proof}
\begin{theorem}
	Let $S$ be a partial semigroup and assume that there is a minimal left ideal $L$ of $S$ which has an idempotent $e$. Then, $e\bullet L$ is a partial group.
\end{theorem}
\begin{proof}
	For every $x\in L$, $L\bullet x$ is left ideal of $S$ and $L\bullet x\subseteq L$ so $L\bullet x=L$ and so there is $y\in L$ such that $y\bullet x=e$. Also, $e\bullet L\bullet e=e\bullet L\subseteq L$, so for every $e\bullet x\bullet e\in e\bullet L\bullet e$ there exists $y\in L(e\bullet x\bullet e)$ such that $y\bullet (e\bullet x\bullet e)=e$, and so $e\bullet (y\bullet (e\bullet x\bullet e))=e$ and $(e\bullet y\bullet e)\bullet (e\bullet x\bullet e)=e$. This implies that for every $x\in e\bullet L\bullet e$ there exists $y\in e\bullet L\bullet e$ such that $y\bullet x=e$. Also, $e$ is the left and right identity for $e\bullet L\bullet e$. Now, by Definition \ref{def-p.g}, $e\bullet L\bullet e$ is a partial group.
\end{proof}
\begin{theorem}
	Let $S$ be a partial semigroup and assume that there is a minimal left ideal $L$ of $S$ which has an idempotent $e$. Then, $L=X\bullet G$ and $X\times G$ are partially isomorphic, where $X$ is the left zero semigroup of idempotents of $L$, and $G=e\bullet L=e\bullet S\bullet e$ is a partial group. 
\end{theorem}
\begin{proof}
	For every $x\in L$, $L\bullet x$ is left ideal of $S$ and $L\bullet x\subseteq L$ so $L\bullet x=L$ and so there is $y\in L$ such that $y\bullet x=e$. By Lemma \ref{2.7}, $e$ is a right identity for $L\bullet e=L$. Therefore by Theorem \ref{1.40} applies with $L$ replacing $S$.
\end{proof}
\begin{theorem}{\label{1.46}}
	Let $S$ be a partial semigroup, let $L$ be a minimal left ideal of $S$, let $T\subseteq S$ and $\bigcap_{s\in L}(R(s)\cap T)\neq\emptyset$. Then $T$ is a minimal left ideal of $S$ if and only if there is some $a\in \bigcap_{s\in L}(R(s)\cap T)$ such that $T= L\bullet a$.
\end{theorem}
\begin{proof}
	Let $T$ be a minimal left ideal. Pick $a\in \bigcap_{s\in L}(R(s)\cap T)$, so $S\bullet (L\bullet a)\subseteq L\bullet a$ and $L\bullet a\subseteq S\bullet T\subseteq T$. Therefore $L\bullet a$ is left ideal of $S$ and $L\bullet a\subseteq T$. This implies that $T=L\bullet a$.
	
	Conversely. Let $T= L\bullet a$ for some $a\in \bigcap_{s\in L}(R(s)\cap T)$. Since $L$ is the left ideal of $S$, so $L\bullet a$ is a left ideal. We show that $L\bullet a$ is a minimal left ideal. So, let $B\subseteq L\bullet a$ be a left ideal of $S$. Define $A=\{s\in L:s\bullet a\in B\}$. Therefore $A\subseteq L$ and $A$ is non-empty. Let $s\in A$ and $t\in L(s\bullet a)$, and so $s\bullet a\in B$ and $t\bullet (s\bullet a)\in B$. Since $s\in L$, $t\bullet s\in L$ and so $t\bullet s\in A$. This implies that $A$ is a left ideal of $S$. Thus $A=L$ and so $L\bullet a\subseteq B$ and $L\bullet a=B$.
\end{proof}
Let $S$ be a partial semigroup. If $S$ has the smallest ideal, we denote  $K(S)$ as the smallest ideal.
\begin{lemma}
	Let $S$ be a partial semigroup, let $I$ be an ideal of $S$, let $L$ be a minimal left ideal of $S$ and let 
	$I\cap L\neq\emptyset$. Then $L\subseteq I$.
\end{lemma}
\begin{proof}
	Pick $x\in L\cap I$. Then $L=S\bullet x\subseteq L\cap I\subseteq I$.
\end{proof}
\begin{theorem}
	Let $S$ be an adequate partial semigroup. If $S$ has a minimal left ideal, then $K(S)$ exists and $K(S)=\bigcup\{L:L\mbox{ is a minimal left ideal of }S\}$.
\end{theorem}
\begin{proof}
	Let $I=\bigcup\{L:L\mbox{ is a minimal left ideal of }S\}$. We show that $I$ is an ideal of $S$. Since there exists a left ideal, $I\neq\emptyset$, let $x\in I$ and $s\in L(x)$. Therefore, for some minimal left ideal $L$, $x\in L$. Then, $s\bullet x\in L\subseteq I$. Now, let $s\in R(x)$. According to Theorem \ref{1.46}, $L\bullet s$ is minimal left ideal of $S$ and $L\bullet s\subseteq I$. This implies that $x\bullet s\in I$, hence $I$ is an ideal.
\end{proof}
\begin{theorem}
	Let $(S, \bullet)$ be a discrete right adequate partial semigroup and let $\mathscr{A} \subseteq \mathcal{P}(S)$ have the finite intersection property. If for each $A \in \mathscr{A}$ and $x \in A,$ there exists $B \in \mathscr{A}$ such that $B \subseteq R(x)$ and $x \bullet B \subseteq A,$ then $\bigcap_{A \in \mathscr{A}} \hat{A}$ is a partial subsemigroup of $\beta S.$
\end{theorem}
\begin{proof}
	Let $T= \bigcap_{A \in \mathscr{A}} \hat{A}.$ Since $\mathscr{A}$ has the finite intersection property, $T \ne \emptyset.$ Let $p, q \in T$, $q \in R(p)$ and let $A \in \mathscr{A}. $ Given $x \in A ,$ there is some $B \in \mathscr{A}$ such that $x \bullet B \subseteq A$ and hence $x^{-1}A \in q.$ Thus $A \subseteq \{x \in S \colon x^{-1}A \in q\}$ so $\{x \in S \colon x^{-1}A \in q\} \in p$, so by Lemma 1.7 in \cite{McLeod1}, $A \in p \bullet q.$
\end{proof}

Let $D$ be a discrete space and $X$ be a compact Hausdorff topological space $X$. Pick $p\in\beta D$ and let $f:D\rightarrow X$ be an arbitrary function. For $y\in X$, we say that  $p-lim_{s\in D}f(s)=y$ if for every neighborhood $U$ of $y$, $f^{-1}(U)\in p$. In fact, if $f^\beta$ is the unique extension of $f$, then $f^\beta(p)=y$. Since $X$ is a compact space, $p-lim_{s\in D}f(s)$ exists and is unique; see Theorem 3.48 in \cite{Hindman}. If $g:X\rightarrow X$ is continuous, then by Theorem 3.49 in \cite{Hindman}, we have 
\[
p-lim_{s\in D}g(f(s))=g(p-lim_{s\in D}f(s)).
\]
Therefore $f^\beta(p)=p-lim_{s\in D}f(s)$ for every $p\in\beta D$.

If $(S,\bullet)$ is a partial semigroup, and $p\bullet q$ is well defined for $p,q\in\beta S$, then 
\[
p\bullet q=p-lim_{s\in P}(q-lim_{t\in Q}s\bullet t),
\]
where $P\in p$ and $Q\in q$.

\begin{lemma}{\label{par-lef-ide}}
	Let $S$ be a discrete partial semigroup, then the following statements hold:
	\begin{flushleft}
		(a) If $L(x)\neq\emptyset$ for some $x\in \beta S$, then  $\overline{L(x)}=\{y\in \beta S:y\bullet x\mbox{ is defined}\}$.\\
		(b) If $L(x)\neq\emptyset$ for some $x\in \beta S$,  $\overline{L(x)}\bullet x$ is left ideal of $\beta S$.
	\end{flushleft}
\end{lemma}
\begin{proof}
	
	(a) Let $p \in \overline{L(x)}$, that is, $L(x) \in p.$ It is easy to see that 
	\[
	p \bullet x = p\mbox{-lim}_{s \in L(x)} s\bullet x.
	\]
	Let us assume that $y = p\mbox{-lim}_{s \in L(x)} s\bullet x$. Then for $A \in y,$ we may choose a neighborhood $U_y: = \overline{A}.$ Now applying the definition of $p$-lim, we have 
	\[
	\{s \in L(x): s\bullet x \in \overline{A}\} \in p.
	\]
	Which implies $\{s \in L(x): s^{-1}A \in x\}\in p$. Since $s \in L(x)$ means $x\in\overline{R(s)}$, therefore, the above containment implies 
	\[
	\{s \in S: x\in\overline{R(s)} \mbox{ and } s^{-1}A \in x\}=\{s \in S : s^{-1}A \in x\}  \in p.
	\]
	This implies $A \in p\bullet x.$ Since we are with ultrafilters, hence $y=p \bullet x.$ This shows that if $p \in \overline{ L(x)}$, then $p \bullet x$ is defined.
	
	(b) Pick $y\in\overline{L(x)}$ and $s\in L(y\bullet x)$, then $y\bullet x$ and $s\bullet (y\bullet x)$ are well defined. Since  $s\bullet (y\bullet x)=(s\bullet y)\bullet x$, we will have $s\bullet y\in \overline{L(x)}$. This implies that $\overline{L(x)}\bullet x$ is left ideal.
\end{proof}
\begin{definition}
	Let $(S,\bullet)$ be a right adequate partial semigroup and let $A\subseteq S$.
	\begin{flushleft}
		(i)
		We say that $A$ is partially syndetic if for every $p\in\beta S$, there is some $H \in \mathcal{P}_f(L(p))$ such that $R(H)\subseteq\bigcup_{t\in H}t^{-1}A$.\\
		(ii)
		We say that $A$ is syndetic if there exists $H \in \mathcal{P}_f(S)$ such that $\delta_R S \subseteq \bigcup_{t\in H}\overline{t^{-1}A}$.\\
		(iii)
		We say that $A$ is partially thick if there exits $p\in\beta S$ such that for every $F\in P_f(L(p))$ there exists $t\in R(F)$ such that $F\bullet t\subseteq A$.\\
		(iv)
		We say that $A$ is $\check{c}$-thick if there exists
		$p \in \delta_R S$ such that $\beta S \bullet p \subseteq \overline{A}$.\\
		(v)
		We say that $A$ is partially piecewise syndetic if for every $s\in S$ there exists $H \in\mathcal{P}_f(R(s))$ such that for every finite nonempty set $T \in R(H)$, there exists $x \in R(T)$ such that $T \bullet x \subseteq \bigcup_{t\in H} t^{-1}A$.\\	
		(vi)
		We say that $A$ is $\check{c}$-piecewise syndetic if there exists  $H \in \mathcal{P}_f(R(s))$  such that for all $T \in \mathcal{P}_f(R(s))$ there exists $x \in R(T)$ such that $(T \cap R(H)) \bullet x \subseteq \bigcup_{t\in H} t^{-1}A$.				
	\end{flushleft}	
\end{definition}

By Lemma 3.5. in \cite{McCutcheon}, we know that $A\subseteq S$ is  $\check{c}$-syndetic if and only if $\overline{A}\cap(\beta S\bullet p)\neq\emptyset$ for every $p\in\delta_R S$. Also Theorem 4.4. in \cite{McCutcheon} say us, $A\subseteq S$ is $\check{c}$-piecewise syndetic if and only if $\overline{A}\cap(\beta S\bullet p)\neq\emptyset$ for some $p\in K(\delta_R S)$.

\begin{theorem}
	Let $(S,\bullet)$ be a right adequate partial semigroup, and let $A\subseteq S$. 
	\begin{flushleft}
		(a) $A$ is partially thick if and only if $S\setminus
		A$ is not partially syndetic.\\
		(b) $A$ is partially syndetic if and only if $S\setminus
		A$ is not partially thick.
	\end{flushleft}
\end{theorem}

\begin{proof}
	(a) Let $S\setminus A$ be partially syndetic. Therefore for every $x\in\beta S$, there exists $H\in P_f(L(x))$ such that $R(H)\subseteq\bigcup_{t\in H}t^{-1}(S\setminus A)$. Since $A$ is partially thick, for some $p\in\beta S$ and for every $H\in P_f(L(p))$ there exists $t\in R(H)$ such that $H\bullet t\subseteq A$. Therefore for $x=p$, we will have  
	\[
	\bigcup_{t\in H}t^{-1}(S\setminus A)\cap \bigcap_{t\in H}t^{-1}A\neq\emptyset.
	\]
	Therefore, $S\setminus A$ can not be partially syndetic.\\
	Conversely, if $A$ is not partially thick, then for every $p\in\beta S$, for some $F\in P_f(L(p))$ and for every $t\in R(F)$ we have $F\bullet t\nsubseteq
	A$, and so $R(F)\cap\bigcup_{t\in F}t^{-1}(A)=\emptyset$. Hence, for every $p\in\beta S$ there exists $F\in P_f(L(p))$ such that $R(F)\subseteq \bigcup_{t\in F}t^{-1}(S\setminus A)$. This implies that $S\setminus A$ is partially syndetic.
	
	(b) It is obvious.
\end{proof}

\begin{theorem}
	Let $(S,\bullet)$ be a right adequate partial semigroup, and let $A$ be a subset of $S$.
	\begin{flushleft}
		(a) $A$ is partially thick if and only if  $\beta S\bullet p\subseteq \overline{A}$ for some $p\in\beta S$.\\
		(b) $A$ is partially syndetic if and only if for every $p\in \beta S$ there exists $H\in P_f(L(p))$ such that  for every $x\in\overline{R(H)}$ there exists $t\in H$ such that $t\bullet x\in\overline{A}$.\\
		(c)  $A$ is partially piecewise syndetic if and only if for every $s\in S$, there exists $G\in P_f(L(s))$ such that $\bigcup_{t \in G}t^{-1}A$ is partially thick.
	\end{flushleft}
\end{theorem}

\begin{proof}
	(a) Since $A$ is partially thick, so $\{t^{-1} A:t\in L(p)\}$ has the finite intersection property for some $p\in\beta S$ and particular, $t^{-1} A\in p$ for every $t\in L(p)$. So $t\bullet p\in\overline{A}$ for every $t\in L(p)$. By Lemma \ref{par-lef-ide}(b), $\overline{L(p)\bullet p}$ is a left ideal, which completes our proof.
	
	Now let $\beta S\bullet p\subseteq \overline{A}$ for some $p\in\beta S$. Therefore $L(p)=\{s\in S:s\bullet p\in\overline{A}\}$. This implies that $\{t^{-1}A:t\in L(p)\}$ has the finite intersection property. So, for every $F\in P_f(L(p))$ there exists $x\in R(F)$ such that $F\bullet x\subseteq A$. This completes the proof. \\
	(b) Let $A\subseteq S$ be a partially syndetic set. So, for every $p\in\beta S$ there exists $H\in P_f(L(p))$ such that $R(H)\subseteq\bigcup_{t\in H}t^{-1}A$. Pick  $x\in\overline{R(H)}$, so for some $t\in H$ implies that $t\bullet x\in \overline{A}$. Therefore  $\overline{A}\cap \beta S\bullet x\neq\emptyset$.  
	
	Conversely. Assume that for every $p\in \beta S$ there exists $H\in P_f(L(p))$ such that for every $x\in\overline{R(H)}$ there exists $t\in H$ such that $t\bullet x\in\overline{A}$. Therefore for every $R(H)\subseteq\bigcup_{t\in H}t^{-1}A$.
\end{proof}
\begin{theorem}\label{4.39}
	Let $S$ be a right adequate partial semigroup and let $p\in\delta_RS$. The following statements are equivalent:
	\begin{flushleft}
		(a) $p\in K(\beta S)$. \\
		(b) For all $A\in p$, $\{s\in S :s^{-1}A\in p\}$ is partially syndetic.\\
		(c) For all $q\in \overline{L(p)}$, $p\in(\beta S\bullet q)\bullet p$.
	\end{flushleft}
\end{theorem}
\begin{proof}
	(a) implies (b). Let $ A \in p$ and let $B=\{s\in S:s^{-1}A\in p\}$. Since $p\in K(\beta S)$, there exists a minimal left ideal $L$ such that $p\in L$. Therefore for every $q\in L$, we have $p \in \overline{S\bullet q}$. Therefore, 
	$s\bullet q\in \overline{A} $ for some $s\in L(q)$, and so $q\in\overline{s^{-1}A}$. Thus, the sets of the form $\overline{s^{-1}A}$ cover the compact set $L$ and hence
	$L\subseteq\bigcup_{s\in H}\overline{s^{-1}A}$ for some finite subset $H$ of $P_f(L(q))$.
	
	Now, let $u\in\beta S$. To see that $R(T)\subseteq\bigcup_{s \in H} s^{-1}B$ for some $T\in P_f(L(u))$, pick $t\in R(T)$, so $t\bullet p\in L$. Then $t\bullet p\in\overline{s^{-1}A}$ for some $s\in H$, and so $s^{-1}A\in t\bullet p$. Therefore, $t^{-1}(s^{-1}A)\in p$ and $s\bullet t\in B$, and thus $t\in s^{-1}B$.
	
	(b) implies (c). Let $q\in \overline{L(p)}$, and suppose that $p\notin(\beta S\bullet q)\bullet p$. Pick $A\in p$ such
	that $\overline{A}\cap(\beta S\bullet q)\bullet p=\emptyset$. Let $B=\{s\in S:s^{-1}A\in p\}$
	and pick $H\in P_f(S)$ such that $S=\bigcup_{s \in H}s^{-1}B$. Pick $s\in H$ such that
	$s^{-1}B \in q $. Then $ B \in s\bullet q $. That is,
	$\{t\in S:t^{-1}A\in p\}\in s\bullet q$ so $A\in(s\bullet q)\bullet p$, a contradiction.\\
	(c) implies (a). Pick $ q \in K(\beta S)$. Therefore
	$p\in(\beta S\bullet q)\bullet p= \beta S\bullet (q\bullet p)$ is minimal left ideal and so $p\in K(\beta S)$.
\end{proof}

\subsection{Right topological left adequate semigroup }
\begin{definition}
	A right topological left adequate partial semigroup is a triple $(S, \bullet, \tau)$, where $(S, \bullet)$ is a left adequate  partial semigroup, $(S, \tau)$ is a topological space and  $r_s:\overline{L(s)}\rightarrow S$ defined by $r_s(x)=x\bullet s$ is continuous for every $s\in S$ and $x\in \overline{L(s)}$. 
\end{definition}
The topological center of a right topological adequate partial semigroup is defined by 
\[
\Lambda(S)=\{s\in S: \lambda_s:\overline{R(s)}\rightarrow S\mbox{ is continuous}\}.
\]
The left topological partial semigroup is defined similarly. 

\begin{theorem}\label{2.3}
	Let $S$ be a compact right topological left adequate partial semigroup,  then $E(S) \ne \emptyset$.
\end{theorem}
\begin{proof}
	Let $\mathscr{A} = \{T \subseteq S \colon  T \ne \emptyset ,~ T ~ \text{is compact, and } ~ \emptyset\neq T \bullet x \subseteq T, ~ \forall x\in T \}$. It means that, $\mathscr{A}$ is the set of all compact partial subsemigroups of $S$.  Since $S \in \mathscr{A}$, so $\mathscr{A} \neq \emptyset$. Let $\mathscr{C}$ be a chain in $\mathscr{A}$. Then $\mathscr{C}$ is a collection of closed subsets of $S$ with the finite intersection property, so  $\bigcap \mathscr{C} \neq \emptyset$ and $\bigcap \mathscr{C}$ is compact. Trivially $\bigcap \mathscr{C} \in \mathscr{A}$. Therefore, $\mathscr{A}$ has a minimal member by Zorn’s Lemma. Let $A \in \mathscr{A}$ be a minimal member. Pick $x \in A$, so  $(A \cap \overline{L(x)})$ is closed and therefore $B=(A \cap \overline{L(x)}) \bullet x \subseteq A \bullet A \subseteq A$ is closed subset of $A$. So $B \bullet B = (A \cap \overline{L(x)}) \bullet x \bullet ((A \cap \overline{L(x)}) \bullet x) \subseteq B$ then $B \in \mathscr{A}$ and $B \subseteq A$ this implies that $B = A$. Therefore, there exists $y \in A$ such that $y \bullet x = x$. 
	Let $C = \{y \in A : y \bullet x = x\}$.
	Let $y, z \in C$. Then $y \bullet z \in A$ and $x=y\bullet x=y\bullet(z\bullet x)=(y\bullet z)\bullet x$. This implies that 
	$y\bullet z\in C$. Thus, $C \in \mathscr{A}$. Since $C \subseteq A$  and $A$ is minimal, it follows that $C = A$. Consequently, 
	$x \in C$ and  $x \bullet x = x$.	
\end{proof}

\begin{theorem}\label{2.4}
	Let $S$ be a compact right topological left adequate partial semigroup. Then, any left ideal of $S$ contains a minimal left ideal. Minimal left ideals are closed, and each minimal
	left ideal has an idempotent.
\end{theorem}
\begin{proof}
	Let $L$ be a left ideal of $S$. Let $\mathscr{A} = \{A \subseteq L : A ~ \text{is a compact left ideal.}\}$. Choose $x \in L$ then $S \bullet x = \rho_x[S]=\overline{ L(x)}\bullet x$ is compact . Also, for every $s\in L(x)$, and for every $t\in L(s\bullet x)$, we have $t \bullet (s \bullet x) = (t \bullet s) \bullet x \in S \bullet x$ so $S \bullet x$ is a left ideal. Since $S \bullet x \subseteq S \bullet L \subseteq L$, thus $\mathscr{A} \ne \emptyset$. We see easily that the intersection of a chain in $\mathscr{A}$ is again in $\mathscr{A}$. So, by Zorn's
	Lemma, has a minimal member $M$ in $\mathscr{A}$, and hence $M $ is a compact minimal left ideal
	contained in $L$. Now, assume that $B \subseteq M$ is a left ideal. Pick $x \in B$, so $S \bullet x =B\in \mathscr{A}$,  while $S \bullet x \subseteq S \bullet B \subseteq B \subseteq M$, so $S \bullet x = M$.  Therefore, $B = M$, and by Theorem \ref{2.3}, every minimal left ideal contains an idempotent.
\end{proof}

\begin{lemma}{\label{leideal}}
	Let $S$ be a compact right topological partial semigroup and let $L(x)\neq\emptyset$ and $R(x)\neq\emptyset$ for every $x\in S$. Then, the following statements hold:
	\begin{flushleft}
		(a) For every $x\in S$, $\overline{L(x)}\bullet x$ is closed left ideal of $S$.\\
		(b) For every $x\in S$, $x\bullet R(x)$ is right ideal of $S$.
	\end{flushleft}
\end{lemma}
\begin{proof}
	By Theorem \ref{par-lef-ide} it is obvious.
\end{proof}
\begin{theorem}{\label{2.20}}
	Let $S$ be a compact right topological partial semigroup, let $L(x)\neq\emptyset$ for every $x\in S$ and let $T$ be a partial subsemigroup contained in the topological center of $S$, then $\overline{T}$ is a partial semigroup.
\end{theorem}
\begin{proof}
	Let $N(x)$ be the collection of all neighborhoods of the point $x$. Let $p,q\in\overline{T}$ and let $q \in R(p)$. For every $U \in N(p \bullet q)$, there exists a neighborhood $V \in N(p)$ such that $r_q(V) \subseteq U$. Therefore, there exists $x_V \in V \cap T$ so that $r_q(x_V) \in U$, and hence $x_V \bullet q \in U$. Thus, there exists $W \in N(q)$ such that $\lambda_{x_V}(W) \subseteq U.$  Consequently, there exists $y_W \in W \cap T$ such that $x_{V} \bullet y_{W} \in U$. This implies that $T \cap U \ne \emptyset$. Finally, $p \bullet q \in \overline{T}.$
\end{proof}

We recall the definition of a partial semigroup homomorphism in the following.
\begin{definition}
	Let $S$ and $T$ be two partial semigroups, and let $f \colon S \to T$ be a function. Then, $f$ is a partial semigroup homomorphism if whenever $x \in S$ and $y \in R_S(x)$, thus for $f(y) \in R_T(f(x))$,  $f(x \bullet y) =f(x) \bullet f(y).$	
\end{definition} 

\begin{theorem}\label{Theorem 4.21 of Hindman}
	Let $(S, \bullet)$ be a right adequate partial semigroup and let $\mathscr{A} \subseteq \mathcal{P}(S)$ have the finite intersection property. Let $(T, \bullet)$ be a compact right topological partial semigroup and let $\varphi \colon S \to T$ satisfy $\varphi(S) \subseteq \Lambda(T).$ Assume that there is some $A \in \mathscr{A}$ such that for each $x \in A$ there exists $B \in \mathscr{A}$ such that $B \subseteq R_S(x),$  $\varphi(x) \in R_T(\varphi(y))$ and $\varphi(x \bullet y) = \varphi(x) \bullet \varphi(y)$ for every $y \in B.$ Then for all $p,q \in \bigcap_{A \in \mathscr{A}} \hat{A},$ if $q \in R(p)$ and $\tilde{\varphi}(p) \in R(\tilde{\varphi}(q)),$ then $\tilde{\varphi}(p \bullet q) = \tilde{\varphi}(p) \bullet \tilde{\varphi}(q)$.
\end{theorem}
\begin{proof}
	Let $p, q \in \bigcap_{A \in \mathscr{A}} \hat{A}$. For each $x \in A,$ Let $q \in \overline{ R(x)}$ and $y \in R(x)$ so, 
	\begin{align*}
	\tilde{\varphi}(x \bullet q) &= \tilde{\varphi}(x \bullet \lim_{y \to q}y) \\
	&= \lim_{y \to q} \varphi(x \bullet y)  ~~~ because ~~ \tilde{\varphi} \circ \lambda_x ~ is ~ continuous\\
	&= \lim_{y \to q} \varphi(x) \bullet \varphi(y) ~~~ because ~ \varphi(x \bullet y)= \varphi(x) \bullet \varphi(y) ~ on ~ a ~ member ~ of ~ q\\
	&= \varphi(x) \bullet \lim_{y \to q} \varphi(y) ~~~ because ~ \varphi(x) \in \Lambda(T)\\
	&= \varphi(x) \bullet \tilde{\varphi}(q) ~ for ~ \tilde{\varphi}(q) \in \overline{R(\varphi(x))}.
	\end{align*}
	Since $A \in p$, let $p \in \overline{L(q)}$,
	\begin{align*}
	\tilde{\varphi}(p \bullet q) &= \tilde{\varphi}((\lim_{x \to p}x) \bullet q) ~for ~ x \in L(q)\cap S\\
	&= \lim_{x \to p} \tilde{\varphi}(x \bullet q) ~ because ~ \tilde{\varphi} \circ \rho_q ~ is ~ continuous\\
	&= \lim_{x \to p} (\varphi(x) \bullet \tilde{\varphi}(q))\\
	&= (\lim_{x \to p} \varphi(x)) \bullet \tilde{\varphi}(q) ~by ~ the ~ continuity ~ of ~ \rho_{\tilde{\varphi}}(q)\\
	&= \tilde{\varphi}(p) \bullet \tilde{\varphi}(q).
	\end{align*}
\end{proof}
\begin{corollary}\label{corollary 4.22 of Hindman}
	Let $(S, \bullet)$ be an adequate partial semigroup and let $\varphi \colon S \to T$ be a partial homomorphism to a compact right topological partial semigroup $(T, \bullet)$ such that $\varphi[S] \subseteq \Lambda(T).$ Then $\tilde{\varphi}$ is a partial homomorphism from $\beta S$ to $T.$
\end{corollary}
\begin{proof}
	By theorem \ref{Theorem 4.21 of Hindman} and for $\mathscr{A}=\{S\}$, the proof is complete.
\end{proof}
Let $(X, \tau)$ be a compact topological space. We define 
\[
\Xi_X^X = \{f : D_f \subseteq X \to X: f \text{ is a function}\}.
\]
For $f, g \in \Xi_X^X $, we define  $g \circ f(x)=g(f(x))$ when $x \in D_f$ and $f(x) \in D_g$. Then $(\Xi_X^X, \circ)$ is a partial semigroup, see \ref{3.1}(1). Also, for every $f \in \Xi_X^X$, we define $L(f)=\{g \in \Xi_X^X \colon R_f \subseteq D_g\}$. In fact, for every $g \in L(f)$, $g \circ f$ is defined. Now, we define $X_{\infty}= X \cup \{\infty\}$. Then $\tau_{\infty}= \tau_X \cup \{\{\infty\}, X \cup \{\infty\}\}$
is a Hausdorff topology on $X_{\infty}$, and  $(X_{\infty}, \tau_{\infty})$ is a compact Hausdorff topological space. It is obvious that, $\infty$ is an isolated point of $X_{\infty}$. Also, for every $f \in \Xi_X^X$, we define $f_{\infty} \colon X_{\infty} \to X_{\infty}$ by 
\begin{align*}
f_{\infty}(x)= \left\{
\begin{array}{rl}
f(x) & x \in D_f\\
\infty & x \notin D_f.
\end{array} \right.
\end{align*}
\begin{lemma}\label{3.1}
	Let $(X, \tau)$ be a compact Hausdorff topological space.
	\begin{flushleft}
		(1)
		$(\Xi_X^X, \circ)$ is a partial semigroup.\\
		(2)
		The function $ \varphi : f \mapsto
		f_{\infty} : \Xi_X^X \to X_\infty^{X_\infty}$
		is bijective.\\
		(3)
		For every $g \in \Xi_X^X$, and every $f \in L(g)$, we have  $\varphi(f \circ g)=\varphi(f) \circ\varphi(g)$.
	\end{flushleft}
\end{lemma}
\begin{proof}
	(1) For $f \in \Xi_X^X$ and $g \in L(f)$, let $D_f = \{x \in X \mid f(x) ~\text{is defined}\}$. If $f \circ (g \circ h)$ is defined then $f \circ(g \circ h)(x) = f(g \circ h(x))= f(g(h(x))) = f \circ g(h(x)) = (f \circ g) \circ h(x)$. This implies that if $f \circ (g \circ h)$ is defined, then $(f \circ g) \circ h$ is defined. Conversely, they should also be defined and then be equal to each other.
	
	(2) Pick $f, g \in \Xi_X^X$ and let $\varphi(f)= \varphi(g)$. Then $f_{\infty}(x)= g_{\infty}(x)$ for every $x \in X_{\infty}$. So, $f(x)=g(x)$ for every $x \in D_f$, this implies that $D_f \subseteq D_g$, and so $D_f = D_g$. This implies that $f=g$ and hence $\varphi$ is one-to-one. Now, we show that $\varphi$ is surjective. Pick $f \in X_\infty^{X_\infty}$, and set $A=f^{-1}(\{\infty\})$. Now, let $g= f \lvert_{A^{c}}$. Then $\varphi(g)=f$, this completes the proof. 
	
	(3) Pick $g \in \Xi_X^X$ and $f \in L(g).$ Then 
	\[
	\varphi(f \circ g)(x) 
	= \left\{
	\begin{array}{cc}
	f(g(x)) & \text{if} ~ g(x) \in D_f\text{ and }x \in D_g\\
	\infty & \text{if } ~ g(x) \notin D_f\text{ or }x \notin D_g.
	\end{array} \right.
	\] 
	Also,
	\[
	\varphi(f) \circ \varphi(g)(x)  = \left\{
	\begin{array}{cc}
	\varphi(f)(g(x)) & \text{if} ~ x \in D_g\\
	\infty & \text{if} ~ x \notin D_g
	\end{array} \right.
	= \left\{
	\begin{array}{cc}
	f(g(x)) & \text{if} ~ g(x) \in D_f\text{ and }x \in D_g\\
	\infty &   \text{if } ~ g(x) \notin D_f\text{ or }x \notin D_g.
	\end{array} \right.	
	\]
	This implies that $\varphi(f \circ g)=\varphi(f) \circ\varphi(g)$.
\end{proof}

Since $\varphi: \Xi_X^X \to {X_\infty}^{X_{\infty}}$ is bijective and $ {X_\infty}^{X_{\infty}}$ is a compact Hausdorff topological space respect to product topology. Let $\tau=\{\varphi^{-1}(U) \colon U~ \text{is open in} ~ X_{\infty}^{X_{\infty}}\}$. Then, $\tau$ is a topology on $\Xi_X^X$, and so $\varphi$ is continuous. Therefore, $\Xi_X^X$ and $X_{\infty}^{X_{\infty}}$ are homeomorphic. We say that a net $\{f_\alpha\}_{\alpha\in I}$ is pointwise convergence to $f$ in $\Xi_X^X$ if for every $x \in \bigcap_{n=1}^{\infty} Dom f_n$
\[
lim_{\alpha}f_\alpha(x)= f(x)\text{ if and only if } lim_{\alpha}\varphi(f_\alpha)(x)= \varphi(f)(x).
\]
\begin{lemma}\label{enveloping right topological semigroup}
	Let $X$	be a compact topological space. Then
	$\Xi_X^X$ is the right topological partial semigroup.
\end{lemma}
\begin{proof}
	Pick $f \in \Xi_X^X$, we show that $r_f \colon L(f) \to \Xi_X^X$ is continuous. Since $r_f(g)= g \circ f$, then $g \circ f(x)= g(f(x))$, for every $x \in D_f$ and $f(x) \in D_g$. Let $g_n \overset{p}{\to} g$ in $\Xi_X^X$, so $\varphi(g_n)(x) \overset{p}{\to} \varphi(g)(x)$ for every $x \in X_{\infty}$. Therefore $\varphi(g_n)(\varphi(f)(x)) \to \varphi(g)(\varphi(f)(x))$ for every $x \in X_{\infty}$, and so $\varphi(g_n) \circ \varphi(f) \overset{p}{\to} \varphi(g) \circ \varphi(f)$. By Lemma \ref{3.1}(3), $\varphi(g_n \circ f) \overset{p}{\to} \varphi(g \circ f)$ in $X_{\infty}^{X_{\infty}}$ and hence $g_n \circ f \overset{p}{\to} g \circ f$ in $\Xi_X^X$.
\end{proof}
\begin{theorem}\label{left translation continuous}
	Let $(X, \tau)$ be a topological space. Then, for each $f \in \Xi_X^X$, $\lambda_f$ is continuous if and only if $f$ is continuous.
\end{theorem}
\begin{proof}
	Suppose that $f$ is continuous. Pick $g, \{g_{\lambda}\}_{\lambda \in I}$ in $L(f)$ such that $g_{\lambda} \overset{p}{\to} g$  in $\Xi_X^X$ then $g_{\lambda}(x) \overset{p}{\to} g(x)$ which is $g(x), g_{\lambda}(x) \in D_f$. Since $f$ is continuous, we have $f(g_{\lambda}(x)) \to f(g(x))$ for every $x \in X$,
	such that $g_{\lambda}(x), g(x) \in D_f.$ Then for $\lambda_f(g_{\lambda}) \overset{p}{\to} \lambda_f(g)$, so $\lambda_f$ is continuous.
	
	Now, let $\lambda_f$ be continuous. Let $\{x_l\}_{l \in I}$ be a net in $D_f$ such that $x_l \to x$ for some $x \in D_f$. We define $g_l(t)=x_l$ and $g(t)=x$ for every $x_l, x \in D_f$ and every $t \in X$. Thus $g_l \overset{p}{\to} g $, by continuity of $\lambda_f$ we have $\lambda_f(g_l) \overset{p}{\to} \lambda_f(g)$. Therefore for every $t \in X$ we can conclude that $\lambda_f(g_l)(t) \to \lambda_f(g)(t)$. Then $f(x_l)= f \circ g_l(t) \to f \circ g(t)=f(x)$, this implies that $f$ is continuous.
\end{proof}

\section{Partial Semigroup Partial Dynamical System}\label{PDS}
\begin{definition}
	A pair $(X,\{T_s\}_{s\in S})$ is called a partial semigroup partial dynamical system (PSPDS) if 
	\begin{flushleft}
		(i)
		$X$ is a compact Hausdorff space, \\
		(ii)
		$(S,\bullet)$ is partial semigroup,\\
		(iii)
		for any $x \in \beta S$, $L(x)=\{s \in S : T_s(x)= s \ast x ~\text{is defined}\}$ is non-empty and the collection $\{L(x):x\in X\}$  has the finite intersection property,\\
		(iv)
		for every $s\in S$, $R(s)=\{x\in X:T_s(x)=s* x\mbox{ is defined}\}=Dom(T_s)$ is a compact non-empty subset of $X$,\\
		(v)
		for every $s\in S$, $T_s:R(s)\rightarrow X$ is continuous, and \\
		(vi)
		for every $s,t\in S$ and $x\in X$, if $s \bullet t$ is defined, then $T_s\circ T_t(x)=T_{s\bullet t}(x)$.
	\end{flushleft}
\end{definition}
\begin{example}
	\begin{flushleft}
		(a) Let $S$ be the generated free semigroup by $A=\{0,1\}$. For any word $u$ over $A$, $|u|$ denotes the number of letters occurring in $u$ and is called the length of $u$. Therefore, $S=\bigcup_{n=1}^\infty A^n$, where $A^n=\{u\in S:|u|=n\}$.  Now, define $\sigma:S\rightarrow S$ by $\sigma(u_1\cdots u_k)=u_2\cdots u_k$ if $k>1$, and is called partially shift function on $S$. Then $(\beta S,\sigma)$ is a partial semigroup partial dynamical syetem. It is obvious that $R(n)=Dom(\sigma^n)=\overline{\bigcup_{k=n+1}^\infty A^k}$ in $\beta S$. \\
		(b) It is obvious that if $S$ is a right adequate partial semigroup, then  $(\beta S, \{\lambda_s\}_{s \in S})$ is a partial semigroup partial dynamical system.
	\end{flushleft}
\end{example}
Let $(X,\{T_s\}_{s\in S})$ be a (PSPDS). For $x\in X$,  $Orb(x)=\{T_s(x): s\in L(x)\}$ is called the orbit of $x$. A non-empty closed subset $Y$ of $X$ is $S$-invariant if $T_s(Y) \subset Y$ for every $s\in S$. If $Y \subseteq X$ is an invariant closed subset of X, we say that $(Y, \{T_s\}_{s \in S})$ is subdynamical system. For any $x \in X$, $\overline{Orb(x)}$ is $S$-invariant. If a  partial semigroup dynamical system $(X, \{T_s\}_{s \in S})$  has no proper subsystems, then we say that it is minimal. A PSPDS $(X,\{T_s\}_{s\in S})$ is minimal if $\overline{Orb(x)}=X$ for every $x \in X$. By Zorn's Lemma, it is obvious that minimal systems exist. 

In the previous definition, it may that $R(s)=X$ for every $s\in S$ in $(X,\{T_s\}_{s\in S})$, but just for a part of $X$, we have 
$T_s\circ T_t=T_{st}$. Therefore, the following definition naturally is defined. 
\begin{definition}
	A pair $(X,\{T_s\}_{s\in S})$ is called a partial semigroup dynamical system (PSDS) if
	\begin{flushleft}
		(a)
		$S$ is a partial semigroup,\\
		(b)
		for every $s \in S$, $R(s)=Dom(T_s)=X$ is a non-empty compact, and \\
		(c)
		for every $s,t \in S$ and $x \in X$, if $(s\bullet t) x$ is defined, then $s (t  x)$ is defined, and thus, they become equal to each other.
	\end{flushleft}
\end{definition}
\begin{example}
	By Example \ref{2.2}, we know that $(P_f(\mathbb{N}),\biguplus)$ is a partial semigroup. Now, let $f:\mathbb{N}\rightarrow S$ be a sequence. Define
	$T_F:S\rightarrow S$ by $T_F(s)=s+\sum_{n\in F}f(n)$, where $F\in P_f(\mathbb{N})$, and we define $T_F\circ  T_G=T_{F\cup G}$ wherever $F\cap G=\emptyset$ for $F,G\in P_f(\mathbb{N} )$. Then $(\beta S,\{T_F\}_{F\in P_f(\mathbb{N})})$ is a partial semigroup dynamical system. 
\end{example}

By Theorem \ref{2.20}, if $(X, \{T_s\}_{s \in S})$ is a PSPDS, then the closure of $\{T_s: s \in S\}$ in $\Xi_X^X$ is a semigroup, which is referred to as the enveloping semigroup of the dynamical system.
\begin{theorem}{\label{teta}}
	Let $S$ be an adequate partial semigroup, $(X, \{T_s\}_{s \in S})$ be a PSPDS and define $\theta : S \to \Xi_X^X$ by $\theta(s)= T_s$. Then $\tilde{\theta}$ is a continuous partial homomorphism  from $\beta S$ onto the enveloping semigroup of $(X, \{T_s\}_{s \in S}).$
\end{theorem}
\begin{proof}
	Since $\tilde{\theta}$ is a continuous extension and $\tilde{\theta}[\beta S]=c\ell \{T_s \colon s \in S\}$. By Lemma \ref{enveloping right topological semigroup}, Theorem \ref{left translation continuous}, each $T_s$ is in the topological center of $\Xi_X^X$.  Therefore, by Corollary \ref{corollary 4.22 of Hindman}, $\theta$ is a partial homomorphism.
\end{proof}
\begin{definition}
	Let $(X, \{T_s\}_{s \in S})$ be a PSPDS and define $\theta : S \to \Xi_X^X$ by $\theta(s)= T_s.$ For each $p \in \beta S,$ let $T_p= \tilde{\theta}(p)$, where $T_p(x)=p-lim_{s \in S} T_s(x).$
\end{definition}
Now, let $(X, (T_s)_{s \in S})$ be a partial semigroup partial dynamical system. Assume that $p\bullet q$ is well-defined for some $p,q\in\beta S$. Pick $x\in X$, and define $f:S\rightarrow X$ by $f(v)=T_v(x)$ for every $v\in S$. Therefore, 
\begin{align*}
(p\bullet q)-lim_{v\in S}f(v)&=p\bullet q-lim_{v\in S} f^\beta(v)\\
&=  f^\beta((p\bullet q)-lim_{v\in S}(v))\\
&=f^\beta(p\bullet q)\\
&= f^\beta (p-lim_{s\in S}q-lim_{t\in S}s\bullet t)\\
&= p-lim_{s\in S}q-lim_{t\in S}f^\beta(s\bullet t)\\
&= p-lim_{s\in S}q-lim_{t\in S}f(s\bullet t).
\end{align*}
Therefore, the following definition is natural. 
\begin{definition}
	Let $(X, (T_s)_{s \in S})$ be a partial semigroup partial dynamical system. For every $x\in X$ and every $p\in \beta S$, we define, 
	\[
	T_p(x)=p-lim_{s\in L(x)}T_s(x),
	\]
	if for every neighborhood $U$ of $T_p(x)$, $\{s \in L(x) : T_s(x)\in U\} \in p$.
\end{definition}
\begin{lemma}\label{4.4}
	Let $q \in \beta S$, $p \in \overline{L(q)}$, then 
	$T_p \circ T_q= T_{p \bullet q}$, where 
	$p \bullet q \in \beta S$ is defined.
\end{lemma}
\begin{proof}
	Pick $p,q\in\beta S$ such that $p\bullet q$ is well defined. Let $x\in R(p\bullet q)$, so 
	\begin{align*}
	T_{p\bullet q}(x) &=(p\bullet q)-lim_{r\in S}T_r(x)\\
	&=q-lim_{t\in S}(p-lim_{s\in S}T_{s\bullet t}(x))\\
	&=q-lim_{t\in S}(T_t(p-lim_{s \in S}T_s(x)))\\
	&=q-lim_{t \in S}(T_t(T_p(x)))\\
	&=T_q(T_p(x)).
	\end{align*}
\end{proof}

\begin{lemma}
	Let $(X, (T_s)_{s \in S})$ be a partial semigroup partial dynamical system. Let $x\in X$ and $p\in\overline{L(x)}$. If $\{s_\alpha\}_{\alpha\in I}$ is a net in $L(x)$ such that $lim_{\alpha}s_\alpha=p$ then $p\ast x=T_p(x)=lim_{\alpha\in I}T_{s_\alpha}(x)$.
\end{lemma}
\begin{proof}
	It is obvious.
\end{proof}

\begin{lemma}
	Let $(X,\{T_s\}_{s\in S})$ be a (PSPDS), then the following statements hold:\\
	(a) For every $x\in X$, $\overline{L(x)}=\{y\in \beta S:y* x\mbox{ is well defined}\}$.\\
	(b) For every $x\in X$, $\overline{L(x)}*x$ is a closed invariant subset of $X$.
\end{lemma}
\begin{proof}
	(a) It is obvious.\\		
	(b) Pick $y\in\overline{L(x)}$ and $s\in L(y*x)$, then $s*(y*x)$ is well defined. Therefore $s* (y* x)=(s\bullet y)* x$, and so $s\bullet y\in \overline{L(x)}$. This implies that $\overline{L(x)}*x$ is a closed invariant subset of $X$.
\end{proof}
\begin{lemma} 
	Let $S$ be an adequate partial semigroup. Then, closed invariant subsets of the partial semigroup partial dynamical system (PSPDS) $(\beta S, \{\lambda_s\}_{s \in S})$ are precisely the closed left ideals of $\beta S$.
\end{lemma}	
\begin{proof}
	Let $I$ be a left ideal of $\beta S$. Pick $p \in I$ and let $s\in L(p)$, then  $s\bullet p\in I$. Therefore, left ideals are invariant. Let $Y$ be a closed invariant subsets of the PSPDS $(\beta S, \{\lambda_s\}_{s \in S})$. Pick $q \in Y$. Since $Y$ is a closed invariant subset of the PSPDS $(\beta S, \{\lambda_s\}_{s \in S})$, we have  
	\[
	O_S(q)=\{s\bullet q:s\in L(q)\}\subseteq Y.
	\]
	This implies that  $Y$ is a left ideal of $\beta S$.
	
\end{proof}
\begin{definition}
	Let $(X, \{T_s\}_{s \in S})$ be a PSPDS. A subset $I$ of $X$ is a minimal closed invariant subset if $I$ is minimal in the set $\{P \subseteq X : P ~\text{is nonempty, closed, invariant}\}.$
\end{definition}

\begin{theorem}
	Let $S$ be an adequate partial semigroup. Then, $L$ is a minimal closed invariant subset of PSPDS $(\beta S, \{\lambda_s\}_{s \in S})$ if and only if $L$  is a minimal left ideal of $(\beta S, \bullet)$.
\end{theorem}
\begin{proof}
	By Theorem \ref{2.4}, it is obvious.	
\end{proof}
\begin{definition}
	Let $S$ be a partial semigroup. We say that $(X, \{T_s\}_{s \in S})$ is a universal minimal PSDS for $S$ if $(X, \{T_s\}_{s \in S})$ is a minimal PSPDS and for every minimal  PSDS $(Y, \{R_s\}_{s \in S})$, there exists an onto continuous function (which is called a dynamical homomorphism) $\varphi : X \to Y $ such that  $R_s \circ \varphi = \varphi \circ T_s$  for each $s \in S$.
\end{definition}

\begin{theorem}
	Let $S$ be an adequate partial semigroup and let $L$ be a minimal left ideal of the PSPDS $(\beta S, \{\lambda_s\}_{s \in S})$. If $(Y, \{R_s\}_{s \in S})$ is a minimal PSPDS and for some $y \in Y$, $L\subseteq\overline{L(y)}$. Then, $(L, \{\lambda_s\}_{s \in S})$ is an extension of the minimal partial semigroup dynamical system $(Y, \{R_s\}_{s \in S})$.
\end{theorem}
\begin{proof}
	Since $(Y, \{R_s\}_{s \in S})$ is a minimal PSPDS, $\overline{L(y)\ast y}=Y$ for every $y \in Y$. Now, define $g:L(y) \rightarrow Y$ by $g(s)= R_s(y)$, and let $\psi=\tilde{g}:\overline{L(y)}\rightarrow Y$ be the unique extension of $g$. Then immediately $\psi$ is a continuous function from $L$ to $Y$.\\			
	Now, let $t\in L(y)$ and $s \in L(t\ast y)$. Then
	\[
	g\circ\lambda_s(t)=g(s \bullet t)= R_{s \bullet t}(y)= R_s\circ R_t(y)= R_s\circ g(t).
	\]
	Now, pick $p\in\overline{L(y)}$ and let $s\in L(p\ast y)$. We show that $\psi(\lambda_s(p))= R_s(\psi(p))$. Therefore 
	\begin{align*}
	\psi(\lambda_s(p))=&  \psi(s \bullet p)\\
	=&\tilde{g}(s \bullet p)\\
	=& \tilde{g}(p-\lim_{t\in S}(s \bullet t))\\
	=& p-\lim_{t\in S}g(s \bullet t)\\
	=& p-\lim_{t\in S}R_{s \bullet t}(y) \\
	= &p-\lim_{t\in S}(R_s \circ R_t (y))\\
	=& R_s(\lim\limits_{t \rightarrow p}R_t(y)) \\
	=& R_s(R_p(y))\\
	=&  R_s(\psi(p)). 
	\end{align*}
	%			for all $s \in L(p)$ and $p \in I.$
	%			\[
	%			\xymatrix@C=1.5cm@R=1.5cm{I \ar[r]^{\lambda_s'} \ar[d]_{\psi} & I\ar[d]^{\psi} \\ Y \ar[r]_{R_s}& Y}
	%			\]
\end{proof}

\section{Partial Enveloping Semigroups}\qquad

Recall $\{0, 1\}^A= \{f \colon A \to \{0, 1\} \colon \mbox{f is a function}\}$ and let $\Omega= \bigcup_{\emptyset \ne A \subseteq S}\{0, 1\}^A$.

\begin{lemma}\label{19.14}
	Let $S$ be an adequate partial semigroup. For each $s \in S,$ define $T_s \colon \{0 ,1\}^{L(s)\bullet s} \to \Omega$ by $T_s(f)(t) =f (t \bullet s)=f \circ \rho_s(t).$ Then $(\Omega, \{T_s\}_{s \in S})$ is a PSPDS. Pick $s \in S$, then for every $p \in \overline{R(s)}$ and $f \in \{0, 1\}^{L(s)\bullet s}$, we have 
	$$(T_p(f))(s)=1 \Leftrightarrow \{t \in S \colon f(t)=1\} \in s \bullet p.$$
\end{lemma}
\begin{proof}
	Let $s \in S$. To demonstrate that $T_s$ is continuous, it suffices to show that $\pi_t \circ T_s$ is continuous for each $t \in L(s).$ So, let $t \in L(s)$ and let $f \in \{0, 1\}^{L(s)\bullet s}.$ Then, 
	$$(\pi_t \circ T_s)(f)= \pi_t(f \circ \rho_s)= (f \circ \rho_s)(t)= f(t \bullet s)= \pi_{t \bullet s}(f).$$ This implies that $\pi_t \circ T_s= \pi_{t \bullet s}(f)$ and, therefore, is continuous.
	Now, let $t \in S$ and $s \in L(t)$. To show that $T_s \circ T_t= T_{s \bullet t}$, let $f \in \{0, 1\}^{L(t\bullet s)}$. Then,
	$$(T_s \circ T_t)(f) = T_s(T_t(f)) = T_s(f \circ \rho_t)=(f \circ \rho_t) \circ \rho_s =f \circ (\rho_t \circ \rho_s)= f \circ \rho_{s \bullet t} =T_{s \bullet t}(f).$$
	
	Finally, let $s \in S$, $p \in \overline{R(s)}$, and $f\in \{0, 1\}^{L(s)}$. Then, 
	\begin{align*}
	T_p(f)(s)=1 & \Leftrightarrow p-lim_{t \in S}T_t(f)(s)=1\\ &\Leftrightarrow p-lim_{t\in S}T_t(f)(s)=1\\ & \Leftrightarrow p-lim_{t \in S}f(s \bullet t)=1 \\
	&\Leftrightarrow \tilde{f}(s \bullet p)=1\\ & \Leftrightarrow \{t \in S: f(t)=1\} \in s \bullet p.
	\end{align*}

\end{proof}
\begin{theorem}
	Let $S$ be an adequate partial semigroup. Then, $\Omega= \bigcup_{\emptyset \ne A \subseteq S}\{0, 1\}^A$ is PSPDS, and for each $s \in S$, define $T_s:\{0 ,1\}^{L(s)\bullet s} \to \Omega$ by $T_s(f)(t)=f\circ \rho_s(t)$. If for every pair of distinct elements $p$ and $q$ of $\overline{ R(s)}$, there is some $s \in S$ such that  $s \bullet p \neq s \bullet q,$ then $\overline{R(s)}$ is topologically and algebrically isomorphic to the enveloping partial semigroup of $(\Omega, \{T_s\}_{s \in S})$.
\end{theorem}
\begin{proof}
	By Lemma \ref{19.14}, $(\Omega, \{T_s\}_{s \in S})$ is a PSPDS. Define $\theta \colon S \to \Omega^{\Omega}$ by $\theta(s)=T_s.$ Then, by Theorem \ref{teta}, $\tilde{\theta}$ is a continuous partial homomorphism from the partial semigroup $\beta S$ onto the enveloping partial semigroup of $(\Omega, \{T_s\}_{s \in S}).$ Thus, we need only show that $\tilde{\theta}$ is one-to-one. So, let $p$ and $q$ be two distinct elements of $\beta S$ and pick $s \in L(p)\cap L(q)$ such that $s \bullet p \neq s \bullet q.$ Pick $A \in s \bullet p \setminus s \bullet q$ and let $\chi_{A}$ be the characteristic function of $A$. Then, by Lemma \ref{19.14}, $T_p(\chi_{A})(s) =1$ and $T_q(\chi_{A})(s)=0$. Therefore, $\tilde{\theta}(p)=T_p \ne T_q = \tilde{\theta}(q).$
\end{proof}
\begin{definition}
	Let $(X, \{T_s\}_{s\in S})$ be a PSPDS.
	\begin{flushleft}
		(a) A point $y \in X$ is uniformly recurrent if and only if for every neighborhood $U$ of $y$, $\{s\in L(y):T_s(y)\in U\}$ is partially syndetic.\\
		(b) Points $x$ and $y$ of $X$ are proximal if and only if there is a net $\langle s_l\rangle_{l \in I}$ in $L(x)\cap L(y)$  such that the nets $\langle T_{s_l}(x) \rangle_{l \in I} $ and $\langle T_{s_l} (y)\rangle_{l \in I}$ converge to the same point of $X.$ 
	\end{flushleft}
\end{definition}
\begin{definition}
	Let $S$ be an adequate partial semigroup. A set $C \subseteq S$ is dynamically central if
	and only if there exist a PSPDS $(X, \{T_s\}_{s \in S})$, points $x,y\in X$, and a neighborhood $U$ of $y$ such that:
	\begin{flushleft}
		(a)
		$y$ is a uniformly recurrent point of $X,$\\
		(b)
		$x$ and $y$ are proximal, and\\
		(c)
		$C=\{s \in S \colon T_s(x) \in U\}$.
	\end{flushleft}
\end{definition}

\begin{lemma}\label{proximal}
	Let $(X, \{T_s\}_{s\in S})$ be a partial semigroup partial dynamical system, and let $x,y \in X$. Then $x$ and $y$ are proximal if and only if there exists a point $p \in \overline{L(x)\cap L(y)}$ such that $T_{p}(x)=T_{p}(y).$
\end{lemma}

\begin{proof}
	Let $z \in X.$ Notice that if $p \in \overline{L(x)\cap L(y)}$, then there is a net $\langle s_l \rangle_{l \in I}$ in $L(x)\cap L(y)$ which converges to $p$, then $T_{p}(z)=\lim_{l \in I}T_{s_l}(z).$ Thus if $p \in \overline{L(x)\cap L(y)}$ and $T_{p}(x)=T_{p}(y),$ then picking any net $\langle s_l \rangle_{l \in I}$ in $L(x)\cap L(y)$ which converges to $p$, one has $\lim_{l \in I} T_{s_l}(x)=\lim_{l \in I} T_{s _l}(y).$
	
	Conversely, assume that $\langle s_l \rangle_{l \in I}$ in $L(x)\cap L(y)$ such that $\lim_{l \in I}T_{s_l}(x)=\lim_{l \in I} T_{s _l}(y).$ Let $p$ be any limit point of $\langle s_l \rangle_{l \in I}$ in $\overline{L(x)\cap L(y)}$. By passing to a subnet we may presume that $\langle s_l \rangle_{l \in I}$ converges to $p.$ Then $T_{p}(x)=\lim_{l \in I} T_{s _l}(x)=\lim_{l \in I} T_{s _l}(y)=T_{p}(y).$ 
\end{proof}

\begin{theorem}\label{19.23}
	Let $(X, \{T_s\}_{s\in S})$ be a partial semigroup partial dynamical system(PSPDS), let $x \in X$ and let $L$ be a minimal left ideal of $\beta S$. The following statements are equivalent:
	
	\noindent (a) The point $x$ is a uniformly recurrent point of $(X,\{T_s\}_{s\in S}).$
	
	\noindent (b) There exists $u \in L$ such that $T_{u}(x)=x.$
	
	\noindent (c) There exists an idempotent $e \in L$ such that $T_{e}(x)=x.$
	
	\noindent (d) There exists $y \in X$ and an idempotent $e \in L$ such that $T_{e}(y)=x$. 
	
\end{theorem}

\begin{proof}
	$(a)\Rightarrow (b)$. Let $v\in L$. Let $\mathcal{N}$ be the set of all neighborhoods of $x$ in $X$. For each $U \in \mathcal{N},$ let 
	\[
	B_U= \{s \in S: T_s(x) \in U\}.
	\]
	Since $x$ is uniformly recurrent, $B_U$ is partially syndetic for every $U\in\mathcal{N}$. So, for some $F_U \in P_{f}(L(v))$, we have  $R(F_U) \subseteq\bigcup_{t\in F_U}t^{-1}B_U $ and $v\in\overline{R(F_U)}$, so there is $t_U \in F_U$ such that $ t^{-1}_{U}B_U \in v$, thus $t_U\bullet v \in \overline{B_U}$. Let 
	$u$ be a limit point of the net $\{t_U\bullet v\}_{U\in\mathcal{N}}$, and notice that $u\in L$. Since $L$ is a left ideal of $\beta S$, so $u \in L$. To see that $T_{u}(x)= x$, we show that $ B_V \in u$ for every $V\in\mathcal{N}$. Suppose that $B_V \notin u$ for some $V\in\mathcal{N}$. Pick $U\subseteq V$ such that $t_U\bullet v\in\overline{S\setminus B_V}$. Hence, $ t_U\bullet v\in\overline{B_U}\subseteq \overline{B_V}$, a contradiction.

	(b)$\Rightarrow$ (c). 
	$L$ is the union of groups, and so $u\in e\bullet L$ for some $e\in E(L)$. Therefore, $T_e(x)=T_e(T_u(x))=T_{e\bullet u}(x)=T_u(x)=x$.
	
	(c) $\Rightarrow$ (d). It is trivial.
	
	(d) $\Rightarrow$ (a). It is obvious that $T_e(x)=T_e(T_e(y))=T_{e\bullet e}(x)=T_e(x)=x$. 
	Let $U$ be a neighborhood of $x$. Pick a neighborhood $V$ of $x$ such that $cl_XV\subseteq U$ and let $A=\{s\in S:T_s(x)\in V\}$. Since $e-lim_{s\in S}T_s(x)=x$, so $A\in e$. Now, by Theorem \ref{4.39}, it implies that $B=\{s\in S:s\bullet e\in \overline{A}\}$ is partially syndetic. Since 
	\[
	T_s(x)=T_s(T_e(x))=T_{s\bullet e}(x)\in\overline{\{T_t(x):t\in A\}}\subseteq \overline{V}\subseteq U
	\]
	for every $s\in B$, we conclude that $\{s\in S:T_s(x)\in U\}$ is partially syndetic.
\end{proof}
\begin{theorem}
	Let $(X, \{T_s\}_{s \in S})$ be a PSPDS and let $x \in X$. Then there is a uniformly recurrent point $y \in\overline{ \{T_s(x) : s \in S\}}$ such that $x$ and $y$ are proximal.
\end{theorem}
\begin{proof}
	Let $L$ be an arbitrary minimal left ideal of $\beta S$ and pick an idempotent $e\in L$. Let $y=T_e(x)$, so $y\in cl\{T_s(x):s\in S\}$. By Theorem \ref{19.23}, $y$ is a uniformly recurrent point of $(X, \{T_s\}_{s \in S})$. Since $T_e(y)=T_{e}(T_{e}(x))= T_{e \bullet e}(x)=T_{e}(x)$, by Theorem \ref{proximal} we conclude that $x$ and $y$ are proximal points.
\end{proof}

\begin{theorem}\label{left ideal}
	Let $(X, \{T_s\}_{s \in S})$ be a PSPDS and let $x, y \in X$. If $x$ and $y$ are proximal, then there is a minimal left ideal $L$ of $\beta S$ such that $T_{u}(x)= T_{u}(y)$ for all $u \in L$.
\end{theorem}
\begin{proof}
	Pick $v \in \overline{L(x)\cap L(y)}$ such that $T_{v}(x)= T_{v}(y)$ and pick a minimal left ideal $L$ of $\beta S$ such that $L \subseteq \beta S \bullet v$. So, let $u \in L$ and choose $w \in L(v)$ such that $u=w\bullet v$. Then, by using Theorem \ref{4.4}, we have 
	\[
	T_{u}(x)=T_{w\bullet v}(x)=T_{w}(T_{v}(x))=T_{w}(T_{v}(y))=T_{w\bullet v}(y)=T_{u}(y).
	\]
\end{proof}

\begin{theorem}\label{19.26}
	Let $(X, \{T_S\}_{s \in S})$ be a PSPDS and let $x, y \in X.$ There is a minimal idempotent $u\in\overline{L(x)}$ such that $T_{u}(x)=y$ if and only if $x$ and $y$ are proximal and $y$ is uniformly recurrent.
\end{theorem}
\begin{proof}
	Since $u$ is minimal idempotent, $L=\beta S\bullet u$ is a minimal left ideal of  $\beta S$, and $u \in L$. Therefore, $T_{u}(x)=y$, and $u \in L$. By Theorem \ref{19.23}, $y$ is a uniformly recurrent point of the PSPDS $(X, \{T_s\}_{s \in S})$. Also,
	\[
	T_{u}(y)=T_{u}(T_{u}(x))= T_{u \bullet u}(x)= T_{u}(x), 
	\]
	and so $x$ and $y$ are proximal points.
	
	Conversely, let $x $ and $y$ be two proximal points, and $y$ is uniformly recurrent point. Then, by Theorem \ref{left ideal}, there is a left ideal $L$ of $\beta S$ such that $T_{u}(x)=T_{u}(y)$ for all $u \in L$. Since $y$ is a uniformly recurrent point, then by Theorem \ref{19.23}, pick an idempotent $e \in L$ such that $T_{e}(y)=y$.		
\end{proof}

We recall the definition of central for adequate partial semigroup in the following.
\begin{definition}
	Let $(S, \bullet)$ be an adequate partial semigroup and let $A \subseteq S.$ Then $A$ is central if and only if there is some minimal idempotent $p$ in $\beta S$ such that $A \in p.$	
\end{definition}
\begin{theorem}\label{19.27}
	Let $S$ be an adequate partial semigroup and let $B \subseteq S.$ Then, $B$ is central if and only if $B$ is dynamically central.
\end{theorem}
\begin{proof}
	Necessity. Let $A$ be an adequate partial semigroup and let $S$ be an adequate partial subsemigroup of $A.$ Let $A = S \cup \{e\}$ where $e$ is a new identity adjoined to $S$ (even if $S$ already has an identity). Let $\Omega= \bigcup_{\emptyset \neq D \subseteq S}\{0, 1\}^D$, for $s \in S$ and  every $t \in L(s)$ define $T_s \colon \{0 ,1\}^{L(s)} \to \Omega$ by $T_s(f)(t) =f (t \bullet s)=f \circ \rho_s(t).$  Then by Lemma \ref{19.14}, $(\Omega, \{T_s\}_{s \in S})$ is a partial semigroup dynamical system. Let $x= \chi_{B},$  the characteristic function of $B \subseteq A.$ Pick a minimal idempotent $u \in \overline{ R(x)}$ such that $B \in u$ and let $y =T_u(x).$ Then by Theorem \ref{19.26}, $y$ is uniformly
	recurrent and $x$ and $y$ are proximal. Now let $U= \{z \in \Omega \colon z(e)=y(e)\}.$ Then $U$  is a neighborhood of $y$ in $\Omega.$ We
	note that $y(e)=1.$ Indeed, $y=T_u(x)$ so $\{s \in S \colon T_s(x) \in U\} \in u$ so choose some $s \in B$ such that $T_s(x) \in U.$ Then $y(e) =T_s(x)(e)=x(e \bullet s)=1.$ Thus given any $s \in S,$   
	\begin{align*}
	s \in B &\Leftrightarrow x(s)=1\\
	& \Leftrightarrow T_s(x)(e)=1\\
	& \Leftrightarrow T_s(x) \in U.
	\end{align*}
	Sufficiency. Choose a PSPDS $(X, \{T_s\}_{s \in S}),$ points $x, y \in X,$ and a neighborhood $U$ of $y$ such that $x$ and $y$ are proximal, $y$ is uniformly recurrent, and $B=\{s \in S \colon T_s(x) \in U\}.$ Choose, by Theorem \ref{19.26}, a minimal idempotent $u \in \overline{ R(x)}$ such that $T_u(x)=y.$ Then, $B \in u.$
\end{proof}

We know that dynamically central sets are closed under supersets, so we state this for adequate partial semigroup in the following.
\begin{corollary}
	Let $S$ be an adequate partial semigroup and let $B \subseteq C \subseteq S.$ If $B$ is dynamically central, then $C$ is dynamically central. 
\end{corollary}
\begin{proof}
	This follows from Theorem \ref{19.27}.
\end{proof}

\end{document}